\documentclass[a4paper,english,10pt,twoside,article]{memoir}
\usepackage[latin1]{inputenc}
\usepackage[T1]{fontenc}
\usepackage{babel}
\usepackage{amsmath,amssymb}
\usepackage{mathtools}
\usepackage[thmmarks]{ntheorem}
\usepackage{empheq}
\usepackage{graphicx}
\makeevenhead{plain}{\thepage}{}{}
\makeoddhead{plain}{}{}{\thepage}
\makeevenfoot{plain}{}{}{}
\makeoddfoot{plain}{}{}{}
\usepackage{lmodern} 
\usepackage[all]{xy} 
\SelectTips{cm}{} 
\usepackage{microtype} 

\makeatletter 
\newcommand*\wt[2][0.1ex]{%
        \begingroup
        \mathchoice{\wt@helper{#1}{#2}{\displaystyle}{\textfont}}
                   {\wt@helper{#1}{#2}{\textstyle}{\textfont}}
                   {\wt@helper{0.53ex}{#2}{\scriptstyle}{\scriptfont}}
                   {\wt@helper{#1}{#2}{\scriptscriptstyle}{\scriptscriptfont}}%
        \endgroup
        #2%
}
\newcommand*\wtb[2][0.1ex]{%
        \begingroup
        \mathchoice{\wt@helper{#1}{#2}{\displaystyle}{\textfont}}
                   {\wt@helper{#1}{#2}{\textstyle}{\textfont}}
                   {\wt@helper{0.34ex}{#2}{\scriptstyle}{\scriptfont}}
                   {\wt@helper{#1}{#2}{\scriptscriptstyle}{\scriptscriptfont}}%
        \endgroup
        #2%
}
\newcommand*\wt@helper[4]{%
        \def\currentfont{\the#41}%
        \def\currentskewchar{\char\the\skewchar\currentfont}%
        \setbox\tw@\hbox{\currentfont#2\currentskewchar}%
        \dimen@ii\wd\tw@
        \setbox\tw@\hbox{\currentfont#2{}\currentskewchar}%
        \advance\dimen@ii-\wd\tw@
        \rlap{\raisebox{-#1}{$\m@th#3\kern\dimen@ii\widetilde{\phantom{#2}}$}}%
}
\makeatother 

\DeclareMathOperator*{\colim}{{colim}}

\renewcommand{\epsilon}{\varepsilon}

\newcommand{\R}{\mathbb{R}}

\renewcommand{\epsilon}{\varepsilon}

\DeclarePairedDelimiter\absv{\lvert}{\rvert}
\DeclarePairedDelimiter\norm{\lVert}{\rVert}

\theoremseparator{.}
\newtheorem*{Theorem}{Theorem}

\newtheorem{Definition}{Definition}[chapter]

\theoremsymbol{$\ast$}
\theorembodyfont{}

\newtheorem{Example}[Definition]{Example}

\theoremsymbol{$\square$}

\newtheorem*{bevis:lem:11}{Proof of Lemma \ref{lem:11}}

\theoremsymbol{``$\square$''}

\begin{document}

\pagestyle{plain}

$\strut$ \\

$\strut$ \\

{\center

 {\Huge \noindent Fibrancy of Symplectic Homology \\
 \noindent In Cotangent Bundles}

$\strut$ \\

 Thomas Kragh\footnote{Supported by Carlsberg-fondet and MIT.}\\
}
$\strut$ \\

$\strut$ \\

\noindent \textbf{Abstract:} We describe how the result in
\cite{mysympfib} extends to prove the existence of a Serre type
spectral sequence converging to the symplectic homology $SH_*(M)$ of
an exact Sub-Liouville domain $M$ in a cotangent bundle
$T^*N$. We will define a notion of a fiber-wise symplectic homology
$SH_*(M,q)$ for each $q\in N$, which will define a graded local
coefficient system on $N$. The spectral sequence will then have page
two isomorphic to homology of $N$ with coefficients in this graded
local system.

\setlength\cftparskip{-15pt}
\tableofcontents

\chapter{Introduction}\label{cha:intro}

Let $N$ be any closed smooth manifold. Let $T^*N \xrightarrow{\pi} N$
denote the cotangent bundle with its canonical projection. We may
define the Liouville 1-form $\lambda_N$ on $T^*N$ by 
\begin{align*}
  \lambda_{N(q,p)}(v) = p(\pi_*(v)), \qquad q\in N, p\in T^*_qN, v\in
  T_{(q,p)}(T^*N). 
\end{align*}
It is well-known that $\omega=d\lambda_N$ is non-degenerate and hence a
symplectic form on $T^*N$.

A (finite type) \textbf{Liouville domain}
$M=(M,\lambda)=(M,\omega,\lambda)$ is a (compact) exact symplectic
manifold ($\omega=d\lambda$) such that the symplectic dual vector
field $X$ of $\lambda$ defined by
\begin{align*}
  \omega(X,-) = \lambda
\end{align*}
points outwards at the boundary $\partial M$. We will only be
concerned with compact Liouville domains, so without further warning
all such will be of finite type. Note that this condition implies that
$\lambda$ pulled back to the boundary $\partial M$ is a contact form.

\begin{Example}
  \label{ex:1}
  Given a Riemannian structure on $N$ the disc bundle
  \begin{align*}
    D_RT^*N =\{ (q,p)\in T^*N \mid \norm{p}\leq R\}
  \end{align*}
  of any radius $R>0$ is a Liouville domain with the 1-form
  $\lambda_N$ defined above.
\end{Example}

An \textbf{Exact Liouville embedding} $(M',\omega',\lambda')
\subset (M,\omega,\lambda)$ is an embedding $M'\subset M$ of
concurrent dimensions, such that $\omega' = \omega_{\mid M'}$ and
there exists a smooth function
\begin{align*}
  f \colon M' \to \R
\end{align*}
such that $\lambda_{\mid M'} - \lambda' = df$.

\begin{Example}
  \label{ex:2}
  Any closed exact Lagrangian sub-manifold $L\subset D_RT^*N$ not
  intersecting the boundary of $D_RT^*N$ has by the Darboux-Weinstein
  theorem an extension to an exact Liouville embedding $D_\epsilon
  T^*L \subset D_RT^*N$ for some small $\epsilon$ and some Riemannian
  structure on $L$.
\end{Example}

In Section~\ref{cha:symp} we will define symplectic homology
$SH_*(M)$ of a Liouville domain $M$. In Section~\ref{cha:fibsymp} We
will describe a fiber-wise version of this. I.e. a symplectic
homology depending on $q\in N$ denoted $SH_*(M,q)$. These will
define a local coefficient system $SH_*(M,\bullet)$ on $N$ and the
following theorem may be thought of as a type of Serre spectral
sequence for symplectic homology.

\begin{Theorem}
  \label{thm:1}
  For any exact Liouville embedding $(M,\lambda) \subset
  (D_RT^*N,\lambda_N)$ there is a spectral sequence
  $(E^r_{*,*},d_r)_{r\geq 1}$ such that
  \begin{itemize}
  \item{it strongly converges to a filtered quotient of $SH_*(M)$ and}
  \item{the second page term $E^2_{n,m}$ is isomorphic to
      $H_n(N;SH_m(M,\bullet))$.}
  \end{itemize}
\end{Theorem}

As will be evident from the construction in Section~\ref{cha:fibsymp}
the local system $SH_*(M,\bullet)$ is trivial if $M \xrightarrow{\pi}
N$ is not surjective. This implies that $SH_*(M)$ is trivial, and thus
generalizes the fact that any exact Lagrangian $L\subset T^*N$
surjects to $N$ (see \cite{MR1090163}). Indeed, in the case of
$D_\epsilon T^*L$ the symplectic homology never vanishes because it is
loop space homology (see \cite{MR1617648}) with possibly twisted
coefficients (see \cite{mysympfib}). In \cite{mysympfib} we also
considered products in the case of $D_\epsilon T^*L$ to get the strong
result that exact Lagrangians are up to a finite covering space lift a
homology equivalence. However, for simplicity we will not consider
these extra structures here.


\chapter{Symplectic Homology of $M\subset T^*N$}\label{cha:symp}

In this section we describe the symplectic homology of a Liouville
domain $M$ exact embedded in $T^*N$. The precise setup used in this
section is important for understanding the construction in
Section~\ref{cha:fibsymp}. One can define symplectic homology
independent of the fact that $M$ is in $T^*N$, but for the purpose of
proving Theorem~\ref{thm:1} the following explicit construction is
convenient. The fact that this actually defines symplectic homology is
known as localization (see e.g. \cite{MR1726235}).

We will need to consider Floer homology (infinite dimensional Morse
homology) of restrictions of the action integral
\begin{align}\label{eq:1}
  A(\gamma)= \int_\gamma \lambda_N - Hdt,
\end{align}
where $H\colon T^*N \to \R$ is a smooth Hamiltonian and $\gamma \colon
I\to T^*N$ is a path. In this section, we will only consider closed
loops $\gamma \colon S^1 \to T^*N$, and we denote the action
restricted to these by $A^{\Lambda}$. It is well-known that the
critical points of $A^\Lambda$ are in 1-1 correspondences with closed
time-1 periodic orbits of the Hamiltonian flow.

The symplectic homology will be defined as a limit of Floer
homologies of $A^\Lambda_s$ associated to a sequence of Hamiltonians
$H^s$. We will define these rather explicitly. So let
\begin{align*}
  W=(1-\epsilon,1] \times \partial M \subset M \subset T^*N
\end{align*}
denote a symplectic collar neighborhood of $\partial M \subset
M$. I.e. the 1-form $\lambda$ is given by $r \lambda_{\mid \partial
  M}$, where $\lambda_{\mid \partial M}$ denotes the pull back of
$\lambda$ to the boundary, and $r$ is the coordinate in
$(1-\epsilon,1]$. We will only need to consider Hamiltonians $H^s$
which are locally constant away from $W$ and which only depend on $r$
on $W$.

Choose any smooth function $f\colon (1-\epsilon,\infty) \to \R$
which is concave such that $f(r)$ tends to $-\infty$ as $r$ tends to
$1-\epsilon$ and $f(r)$ is $0$ for $r\geq 1$. We then define a smooth
family of smooth functions $f_s\colon \R_{\geq 0} \to \R$ for $s>1$ such 
that
\begin{itemize}
\item {$f_s(r) = f(r)+s$ when $f(r)\geq -s/2$,}
\item {$f_s''$ has a unique 0 in $(1-\epsilon,1)$, and}
\item {$f_s(r)=0$ when $0\leq r <(1-\epsilon)$.}
\end{itemize}
Notice that these imply that $f_s$ is convex on the interval from
\begin{figure}[ht]
  \centering
  \includegraphics{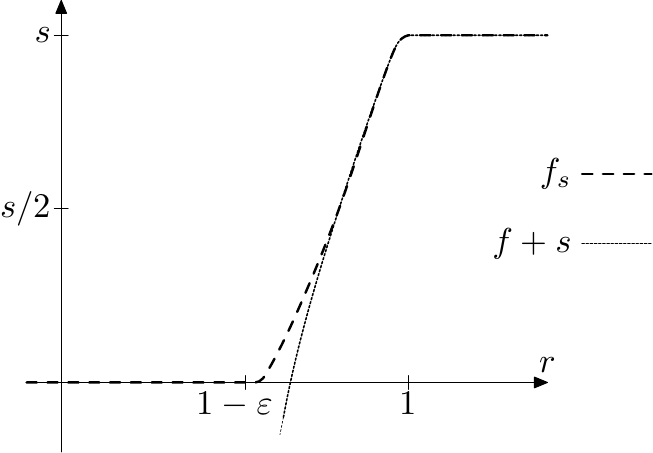}
  \caption{Functions $f_s$ and $f+s$.}
  \label{figfs}
\end{figure}
$-\infty$ to the unique $0$ (in the second bullet point) and concave
on the complement. As illustrated in Figure~\ref{figfs} this describes
a controlled way to ``cap of'' $f+s$ making it smooth, bounded,
non-negative, and more we will need later. We then define the smooth
family of smooth Hamiltonians $H^s$ by
\begin{align*}
  H^s(z) = \left\{
    \begin{array}{ll}
      f_s(r) & z=(r,x) \in W \\
      0      & z\in M-W \\
      s      & z\in T^*N-M
    \end{array}
  \right.
\end{align*}
By construction the Hamiltonian flow of $H^s$ has only constant
periodic orbits outside of $W$ and hence the action on closed loops
not in $W$ only has critical values $0$ and $-s$.

There is a geometric interpretation of the action of other
periodic orbits. Indeed, since we only consider loops in this section
and the embedding $M\subset T^*N$ is exact we may as well integrate
over $\lambda$ and not $\lambda_N$ in Equation~\eqref{eq:1}. Also, if
$\gamma \colon S^1 \to \R$ is a Hamiltonian flow curve it has to have
constant $r$-factor. Indeed, the Hamiltonian is preserved under
Hamiltonian flow. Considering variations of this $r$-factor (knowing
that $\gamma$ is a critical point for $A^\Lambda$) yields the relation
\begin{align} \label{eq:2}
  \int_\gamma \lambda_{\mid \partial M} = f_s'(r),
\end{align}
and thus one calculates that
\begin{align*}
  A^\Lambda_s(\gamma) = rf_s'(r) - f_s(r) = -(f_s(r) - rf'_s(r)),
\end{align*}
which is minus the intersection of the 2. axis with the tangent of
$f_s$ at the point $(r,f_s(r))$.

Now pick $s_0>>1$ large enough such that; the tangent of $f+s_0$ at
the point where $f+s_0 = s_0/2$ intersects the 2. axis below
zero. This and the assumptions on $f_s$ implies that for
$s>s_0$ we thus have; any tangent of $f_s$ intersecting the 2. axis
above $0$ must be a tangent on the part where $f_s=f+s$ (and
$f>-s/2$).

This means that any 1-periodic orbit with negative action must lie on
an $r$-level for which $f_s(r)=f(r)+s$ (and also for $r$'s close to
it). So when increasing $s$ we simply increase $f_s$ by the same value
in a neighborhood. This implies that the 1-periodic orbit is unchanged
and the action is simply translated down with the same speed as $s$ is
translated up. Thus the critical set of $A^\Lambda_s$ has its critical
values \emph{below} zero pushed downwards with the same speed as $s$
increases.

This means that for any $a<0$ we in fact have continuation maps to
higher and higher $s$ and may in fact define the limit
\begin{align}\label{eq:4}
  SH_*(M) = \colim_{s\to \infty} FH_*^a(A_s^\Lambda),
\end{align}
where $FH_*^a(A_s^\Lambda)$ denotes Floer homology associated to
critical values of $A^\Lambda_s$ with critical value greater than
$a$. Indeed, sliding a critical value below $a$ corresponds to
collapsing a generator in the Floer chain complex, which is always a
chain map and which is why the continuation maps exist. Defining this
colimit in the standard way requires picking a sequence of $s$'s
tending to infinity and perturbing the associated actions (and some
compactness proofs). However, to prove Theorem~\ref{thm:1} we will not
use standard Floer theory. We will 
instead use finite dimensional approximations. The advantages of this
is addressed in Section~\ref{cha:proof}. Note that the definition of
this limit does not depend on $a<0$ nor the choice of sequence $s$
going to infinity. Indeed, all critical points with critical values in
$]-\infty,0[$ for some $s$ will eventually get their critical values 
translated down below any given value $a<0$.


\chapter{Fiber-wise Symplectic Homology}\label{cha:fibsymp}

In this section we elaborate and extend the calculation of the action
of 1-periodic orbits to open flow curves. We then use this to
define a fiber-wise version of the symplectic homology.

Let $\gamma \colon I \to T^*N$ be any time-1 flow path for
the Hamiltonian flow of $H^s$. The argument proving
Equation~\eqref{eq:2} can be extended to this case using bump
functions in the variable $u$. I.e. the equation is still
valid. However, when using this to calculate $A_s$ we get a
correction term due to the fact that integrating $\lambda$
and $\lambda_N$ on open paths does not give the same result. Indeed,
let $f \colon M\to \R$ be such that $df= \lambda_N - \lambda$ then
\begin{align*}
  A_s(\gamma) = \int_\gamma \lambda_N -H^sdet = \int_\gamma \lambda
  -H^sdt + f(\gamma(1))-f(\gamma(0)).
\end{align*}
We conclude that the geometric interpretation using intersections of
tangents of $f_s$ with the 2. axis is still valid except we have to
add the term $f(\gamma(1))-f(\gamma(0))$.

Now fix a $q\in N$ and look at the action integrals $A_s^q$ defined on
paths $\gamma$ starting and ending on the Lagrangian fiber
$T_q^*N$. It is well-known that the critical points of $A_s^q$ are
time-1 Hamiltonian flow curves starting and ending on $T^*_qN$.
This means that we can in fact conclude that if $a<\max_{z,z' \in M}
\absv{f(z)-f(z')}$ then the critical set of $A_s^q$ below $a$ slides
down when $s$ is increased. Indeed, as in the previous section the
time-1 flow curves do not change but their action is translated along
with $s$. So we may define
\begin{align}\label{eq:3}
  SH_*(M,q) = \colim_{s\to \infty} FH_*^a(A_s^q).
\end{align}
This is what we will consider the fiber-wise symplectic homology of
$M$. By definition we see that if $M \to N$ does not surject there
must be a fiber on which $H^s$ is constantly equal to $s$. This
implies that for this fiber the fiber-wise symplectic homology is
trivial, and thus by theorem 1 all of them are trivial including
the global symplectic homology $SH_*(M)$.


\chapter{Fibrancy and Sketch of Proof of
  Theorem~\ref{thm:1}}\label{cha:proof}

The proof of Theorem~\ref{thm:1} relies on the fact that the family of
fiber-wise symplectic homologies actual looks and behaves like a Serre
fibration with total space giving the total symplectic homology
$SH_*(M)$. To prove this using homological algebra and infinite
dimensional Floer theory seem very delicate and complicated -
especially in light of transversality issues. So this is not the
approach taken in \cite{mysympfib}. Indeed, there we use finite
dimensional approximations.

\begin{figure}[ht]
  \centering
  \includegraphics{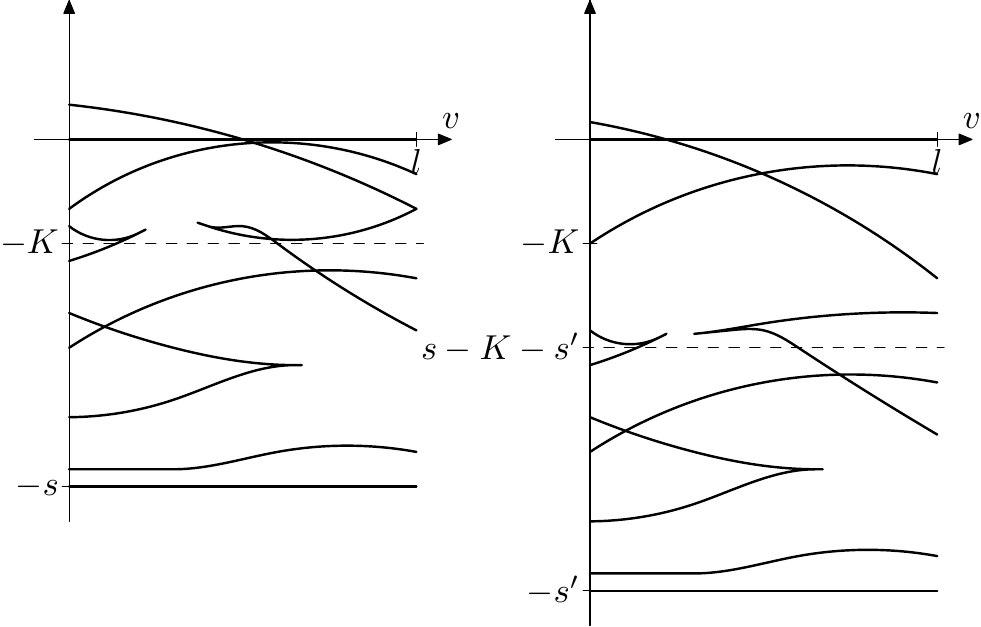}
  \caption{Bifurcation diagrams for same path but different $s$.}
  \label{bifur}
\end{figure}

However, to see why the fiber-wise symplectic homologies form a local
system and understand the construction it is convenient to discuss why
a smooth path $\gamma \colon [0,l] \to N$ (assumed to be parametrized
by arc-length) can induce a ``parallel transport'' in the fiber-wise
homologies. So in the following we describe why such a path gives rise
to a map between the fiber-wise symplectic homologies
\begin{align}\label{eq:5}
  SH_*(M,\gamma(0)) \to SH_*(M,\gamma(l)).
\end{align}
First consider any $s>>1$ fixed. Then for each $v\in [0,l]$ we have
the fiber-wise action functional $A_s^{\gamma(v)}$ defined using the fiber
$T^*_{\gamma(v)}N$ and one may consider the ``graph'' or
\textbf{bifurcation diagram} of the critical set of each as a
multi-valued function of $v$. Examples are illustrated in
Figure~\ref{bifur}.
We know from the discussion in the previous section that the
dependence on $s$ of this diagram is such that everything below some
$-K$ is simply translated downwards with the same speed as $s$ is
translated upwards. Figure~\ref{bifur} shows two possible bifurcation
diagrams for the same system and path $\gamma$, but for different
$s$'s. The part above $-K$ can behave arbitrarily as $s$ changes
except there is always an upper bound. This made the colimits of Floer 
homologies in the previous sections well-defined.

The main point is that the slopes of the ``graph'' pieces in the
bifurcation diagrams are bounded by a certain number. Indeed, Lemma
9.2 in \cite{mysympfib} gives this bound to be 2 if $M\subset DT^*N$,
which can always be arranged by scaling $M$ at the very
beginning. Note that the proof there depends on having a Hamiltonian
which have a slight slope at infinity, but in this heuristical
argument we don't care about the constant critical points outside of
$M$ since we simply assume $s$ very large so that they do not
contribute to the Floer homology (i.e. the critical value $-s$ is much
smaller than any $a<-K$ chosen as the energy cut-off).

The effect of such a bound is that if we make $s$ actually depend on
$v$ and grow faster than $2v$, say $s(v)=s_0+3v$, then the bifurcation
diagram depending on $v$ becomes slanted. Indeed, the slopes of any
critical values with value less than $-K$ has to lie in the interval
$[-5,-1]=[-2,2]-3$. The upshot is that whenever a critical value of
$A_{s(v)}^{\gamma(v)}$ co-insides with $a$ it must go down (as a
function of $v$) with positive speed. We are thus again collapsing
generators in the Floer chain complex and have maps similar to the
continuation maps described in the previous sections. I.e. we
get maps
\begin{align*}
  FH_*^a(A_{s_0}^{\gamma(0)}) \to FH_*^a(A_{s_0+3l}^{\gamma(l)}),
\end{align*}
which when taking the limit as $s_0$ tends to infinity defines the
wanted map from Equation~\eqref{eq:5}. Usual zig-zagging arguments
can be employed to prove that this is an isomorphism, but there are of
course many things to check here, and proceeding using standard Floer
theory seems very cumbersome.

The proof of Theorem~\ref{thm:1} is carried out in \cite{mysympfib}
using the theory of Conley indices to produce spaces instead of
homology theories. I.e. the spaces has the wanted homology. This makes
it possible to define these fiber-wise symplectic
homologies even when $a$ from Equation~\eqref{eq:3} is not regular for
a particular fiber $T^*_qN$ and some $s$. Indeed, this is needed
because there is no guarantee that we can find a sequence of $s$'s
tending to infinity and being regular for all fibers
simultaneously. In fact, in many cases this can be proven not to
exist.

Working with spaces also has the advantage that it is possible to
prove a Serre type fibration property in a more conventional manor,
and thus lifting more than simply a path parametrized by arc length,
but any compact family of smooth paths. It is important here to
mention a technical but significant complication; when increasing $s$
these finite dimensional approximations 
gets more and more complicated and suspensions of the Conley indices
are introduced. This effectively means that the colimits analogous to
those in Equation~\eqref{eq:4} and Equation~\eqref{eq:3} are not
taken in the category of spaces, but in the category of spectra - and
in the fiber-wise case even in a category of parametrized spectra over
$N$. However, the notions of Serre fibrancy, taking fiber-wise
homology, and taking global homology still exist. So we, indeed, get a
Serre type spectral sequence as in Theorem~\ref{thm:1}. We should
note, however, that in \cite{mysympfib} only the case of $M=DT^*L$ is
considered, but all the ideas generalize, and using the Hamiltonians
defined here all the methods to get the spectral sequence immediately
apply. However, to incorporate the products, which is very important
in \cite{mysympfib} one needs to construct the Hamiltonian family
$H^s$ more carefully. 


\bibliographystyle{plain}
\bibliography{../../Mybib}

\end{document}